\documentclass[12pt,reqno]{article}
\usepackage{amssymb}
\usepackage{amsmath}
\usepackage{amsfonts}
\usepackage{amsthm}

\newtheorem{definition}{Definition}
\newtheorem{note}{Note}
\newtheorem{statement}{Statement}
\newtheorem{theorem}{Theorem}

\newtheorem{example}{Example}

\begin{document}
\title{Some Extensions to Touchard's Theorem on Odd Perfect Numbers}
\author{Paolo Starni}
\date{}

\maketitle

\begin{abstract}
The multiplicative structure of an odd perfect number $n$, if any, is $n=\pi^\alpha M^2$, where $\pi$ is prime, $\gcd(\pi,M)=1$ and  $\pi\equiv \alpha\equiv1\pmod{4}$. An additive structure of $n$, established by Touchard, is that  ``$\bigl(n\equiv 9\pmod{36}\bigr )$ OR $\bigl (n\equiv1\pmod{12}\bigr )$''. A first extension of Touchard's result is that the proposition ``$\bigl(n\equiv x^2\pmod{4 x^2}\bigr )$ OR $\bigl (n\equiv \pi\equiv1\pmod{4 x}\bigr )$'' holds for $x=3$ (the extension is due to the fact that the second congruence contains also $\pi$). We further extend the proof to $x=\alpha+2$, $\alpha+2$ prime, with the restriction that the congruence modulo $4 x$ does not include $n$. Besides, we note that the first extension of  Touchard's result holds also with an exclusive disjunction, so that $\pi\equiv 1\pmod{12}$ is a sufficient condition because $3\nmid n$.
\end{abstract}
 
\section{Introduction}
Without explicit definitions all the numbers considered here must be taken as strictly positive integers. 
\begin{definition}
$n$ is said to be perfect if and only if $\sigma(n)=2n$, where $\sigma(n)$ is the sum of the divisors of $n$.
\end{definition}
Euler \cite[p. 19]{Dickson} established a multiplicative structure of odd perfect numbers: 
\begin{statement}[Euler]
If $n$ is an odd perfect number, then $n=\pi^\alpha M^2$, where
$\pi$ is prime, $\gcd(\pi,M)=1$ and $\pi\equiv \alpha\equiv1\pmod{4}$. 
\label{st1}
\end{statement}

In what follows you have to consider that $n$, with the notation used in Statement \ref{st1}, is an odd perfect number and that the proofs of all the theorems are given in Section 2.
\begin{note}[factors of $n$ related to $\sigma(\pi^{\alpha})$]
We can find factors of $n$ considering that:
\begin{equation*}
n=\pi^{\alpha}M^2=\frac{\sigma(\pi^{\alpha})}{2}\sigma(M^2)
\end{equation*}
and
\begin{equation}
\sigma(\pi^{\alpha})=(\pi+1)(1+\pi^2+\pi^4+...+\pi^{\alpha -1}).
\label{eqn:1}
\end{equation}
Since $\gcd(\pi^{\alpha}, \sigma(\pi^{\alpha}))=1$, we have that $\frac{\pi+1}{2}\mid M^2$. In particular, if $\frac{\pi+1}{2}$ is\vspace{0.1cm} squarefree, then $(\frac{\pi+1}{2})^2\mid M^2$.
\label{n1}
\end{note}

With regard to the additive structure of $n$, related to the division by $4$, it holds $n\equiv 1\pmod{4}$; in fact, in Statement 1, $M^2\equiv 1\pmod{4}$ because square of an odd integer.

Touchard  \cite{Touchard} found the additive structure of $n$ in relation to the division by 3: 
\begin{statement}[Touchard]
\begin{equation*}
\bigl (n\equiv 9\pmod{36}\bigr) \hspace{0.3cm} OR \hspace{0.3cm}\bigl ( n\equiv1\pmod{12}\bigr). 
\end{equation*} 
\label{st2}
\end{statement}
Statement \ref{st2} may be rewritten in a little more extended form:
\begin{theorem}
$\bigl (n\equiv 9\pmod{36}\bigr ) \hspace{0.3cm}$ OR $\hspace{0.3cm}\bigl ( n\equiv\pi\equiv 1\pmod{12}\bigr )$.
\label{th1}
\end{theorem} 
Given $\alpha$ as in Statement \ref{st1}, we obtain:
\begin{theorem}
If $\alpha+2$ is prime, then
\begin{equation*} 
\bigl (n\equiv(\alpha+2)^2\pmod {4(\alpha+2)^2}\bigr )\hspace{0.3cm} OR \hspace{0.3cm} \bigl (\pi\equiv1\pmod{4(\alpha+2}\bigr ).
\end{equation*}
\label{th2}
\end{theorem}
 Theorem \ref{th2} reproduces, for $\alpha=1$, Theorem \ref{th1}, except for the fact that its second congruence does not contain $n$.
\begin{note}[factors of $n$ related to $\pi^{\alpha}$]
We can research factors of $n$, in addition to the method of the Note \ref{n1}, considering Theorem \ref{th2}:

 if $\pi\not\equiv 1\pmod{4(\alpha+2)}$ and $\alpha+2$ is prime, then it holds the first congruence of the thesis, so that $(\alpha+2)^2\mid n$. 
\label{n2}
\end{note}
Finally, the most extended form of Toucard's theorem is the alternative:
\begin{theorem}
$\bigl (n\equiv 9\pmod{36}\bigr ) \hspace{0.3cm}$ XOR $\hspace{0.3cm}\bigl ( n\equiv\pi\equiv 1\pmod{12}\bigr )$.
\label{th3}
\end{theorem}
Theorem \ref{th3} has a non-trivial corollary:
\begin{theorem}
If $\pi\equiv 1\pmod{12}$, then $3\nmid n$.
\label{th4}
\end{theorem}
 We apply the previous considerations to the following two examples.
 \begin{example}
 Let it be $n=13^9M^2$. Thus, $7^2\mid n$ from Note \ref{n1} , $11^2\mid n$ from Note \ref{n2} and $3 \nmid n$ from Theorem \ref{th4}.
 \end{example}
 \begin{example}
  Let it be $n=5M^2$. The methods described in Notes \ref{n1} and \ref{n2} lead to the same result that $9\mid n$. 
 \end{example}

\section{The proofs}
We warn that a line numbered to the right, regardless of content, is generically indicated as equation.
\subsection{Proof of Theorem \ref{th1}}
\begin{proof}
Also synthesizing the simpler proofs of Touchard's theorem given by Satyanarayana \cite{Satyanarayana}, Raghavachari \cite{Raghavachari} and Holdener \cite{Holdener}, we prove the following two statements contained in Equations \eqref{eqn:2} and \eqref{eqn:3}:
\begin{equation}
if \hspace{0.1cm} 3\mid n, then\hspace{0.1cm} n\equiv 9\pmod{36}
\label{eqn:2}
\end{equation}
In fact, $3\nmid \pi$ so $9\mid M^2\implies n=\pi^{\alpha}9N^2=(4k+1)9\equiv 9\pmod{36}$.
\begin{equation}
if\hspace{0.1cm} 3\nmid n, then\hspace{0.1cm} n\equiv \pi\equiv 1\pmod{12}
\label{eqn:3}
\end{equation}
In fact, $\pi\neq 12k+5$ (otherwise $3\mid n$, see Equation \eqref{eqn:1}), so $\pi\equiv1\pmod{12}$. Since $M\equiv 1\pmod{6}$ or $M\equiv 5\pmod{6}$, it follows $M^2\equiv 1\pmod{6}$. Thus, being also $M^2\equiv 1\pmod{4}$, it results $ M^2\equiv1\pmod{12}$ and, therefore, $ n\equiv1\pmod{12}$.

 Combining the statements in Equations \eqref{eqn:2} and \eqref{eqn:3}, it follows the proof.
 \end{proof}

\subsection{Proof of Theorem \ref{th2}}

 We need the following result due to Starni \cite{Starni2}:
\begin{statement}
If $\gcd(\pi-1,\alpha+2)=1$ and $\alpha+2$ is prime, then $(\alpha+2)\mid M^2$. 
\label{st4}
\end{statement}
Now we are able to prove Theorem \ref{th2}.
\begin{proof} There are two cases:
\begin{description}
\item $case\hspace{0.1cm} 1)\hspace{0.2cm}\gcd(\pi-1,\alpha+2)=1$.
\item We obtain from Statement \ref{st4}:
\begin{equation*}
n=\pi^{\alpha}(\alpha+2)^2N^2=(\alpha+2)^2(4k+1). 
\end{equation*}
\item It follows the first congruence of the thesis.
\item $case\hspace{0.1cm} 2)\hspace{0.2cm} \gcd(\pi-1,\alpha+2)>1$, i.e., $\gcd(\pi-1,\alpha+2)=\alpha+2$. 
\item We have 
\begin{equation*}
\pi-1=k(\alpha+2)\implies\pi\equiv1\pmod{\alpha+2}.
\end{equation*} 
\item Since $\pi\equiv1\pmod{4}$, it follows the second congruence of the thesis.
\end{description}
\end{proof}

\subsection{Proof of Theorem \ref{th3}}
\begin{proof}
The system of the two congruences in Theorem \ref{th1} does not have solution because $\gcd(36,12)=12\nmid(9-1)$. Thus, the logical connective OR may be replaced by XOR.
\end{proof}

\subsection{Proof of Theorem \ref{th4}}
\begin{proof}
Theorem \ref{th3} states that:
\begin{equation*}
\pi\equiv 1\pmod{12}\implies n \not\equiv  9\pmod{36} 
\end{equation*}
The contrapositive formulation of the statement in Equation \eqref{eqn:2} is
\begin{equation*}
 n \not\equiv  9\pmod{36}\implies 3\nmid n
 \end{equation*}
 Thus, from the transitive property of the logical implication, it follows the proof.
\end{proof}


\begin{thebibliography}{99}

\bibitem{Dickson}
L. E. Dickson, \emph{History of the theory of numbers}, vol. 1, Dover, New York, (2005).
\bibitem{Holdener}
J. A. Holdener, A theorem of Touchard on the form of odd perfect numbers, \emph{ American Mathematical Monthly} \textbf{109} (7) (2002), 661-663.
\bibitem{Raghavachari}
M. Raghavachari, On the form of odd perfect numbers, \emph{ Math. Student} \textbf{34} (1966), 85-89.
\bibitem{Satyanarayana}
M. Satyanarayana, Odd perfect numbers, \emph{Math. Student} \textbf{27} (1959), 17-18.
\bibitem{Starni2}
P. Starni, Odd perfect numbers: a divisor related to the Euler's factor, \emph{J. Number Theory} \textbf{44} (1993), 58-59.
\bibitem{Touchard}
J. Touchard, On prime numbers and perfect numbers, \emph{Scripta Mathematica} \textbf{19} (1953), 35-39.
\end{thebibliography}
\end{document}